\documentclass{amsart}

\usepackage{enumerate}

\usepackage{amssymb}
\usepackage{amsfonts}
\usepackage{amsmath}

\makeatletter
\@namedef{subjclassname@2010}{%
  \textup{2020} Mathematics Subject Classification}
\makeatother

\newtheorem{theorem}{Theorem}[section]
\newtheorem{lemma}{Lemma}[section]

\newtheorem{corollary}{Corollary}[section]

\theoremstyle{definition}
\newtheorem{definition}{Definition}[section]
\newtheorem{example}{Example}[section]
\newtheorem{remark}{Remark}[section]
\newtheorem{proposition}{Proposition}[section]

\usepackage{tikz}

\usetikzlibrary{calc}

\def\afrac#1#2{\ifinner {#1}/{#2} \else \frac{#1}{#2} \fi}
\def\1#1{\mbox{\rm{#1}}}

\allowdisplaybreaks

\begin{document}
\title[The enumeration of plane (3,6)-triangulations] 
 {The enumeration of plane (3,6)-triangulations  and the form of odd perfect numbers}
\author[Jan Florek]{Jan Florek}

\address{Jan Florek, Faculty of Pure and Applied Mathematics,
Wroclaw University of Science and Technology, Wybrze\.{z}e Wyspia\'nskiego 27 \\  50--370 Wroc{\l}aw, Poland}
\email{jan.florek@pwr.edu.pl}

\maketitle
\date{}

\begin{abstract} 


 Let $\mathcal{P}$ be the family of all $3$-connected plane triangulations with vertices of degree $3$ or $6$. The number $d(n)$ of non-isomorphic triangulations belonging to $\mathcal{P}$ of degree $2n+2$ are enumerated,  for $n \geqslant 1$. From the formula for $d(n)$ (where $n$ is odd and $n > 1$), we obtain a strengthening of Euler's theorem regarding the form of odd perfect number.

 \end{abstract}
\bigskip


\keywords{plane triangulation, orbits, enumeration, Euler theorem on odd perfect numbers.

\maketitle

\section{ Introduction}
\bigskip

A connected plane graph $G$ is called a triangulation if every face of $G$ 
(including the outer face) is bounded by a triangle. Let $e$ be an edge of a triangulation $G$. There are exactly two triangles containing $e$. The two vertices of these triangles which do not belong to~$e$ are said to be opposite w.r.t. edge $e$. If a $4$-colouring $c$ (proper) of a triangulation~$G$ is such that any two opposite vertices have different colours, then $c$ is said to be  \textit{nonsingular}. A triangulation has at most one nonsingular colouring. The following was proved by Fisk~\cite{flo3}. 
\begin{proposition}\label{prop1.1} [Fisk]
A plane triangulation has a nonsingular $4$-colouring if and only if the degree of each vertex is divisible by three.
\end{proposition}
Let $\mathcal{P}$ be the family of all simple plane triangulations with  all vertices of degree $3$ or~$6$. Every simple triangulation with at least four vertices is $3$-connected (see Diestel \cite[Corollary 4.4.7]{flo2}). Notice that any two non-simple plane triangulations of the same order with all vertices of degree $3$ or $6$ are isomorphic (see condition (2) of Proposition \ref{prop4.1}).
 
 Let $P \in \mathcal{P}$. Since $P$ is a plane triangulation, then for every vertex $w \in V(P)$ there exists a cyclic orientation around $w$ of all edges which are incident with $w$. Let $g_{0}$, $g_{1}$, $g_{2}$  be fixed edges in $P$  (indexed by elements of the cyclic group $Z_3$) having  counter-clockwise orientation around the common vertex (say $v$) of degree~$3$. 
Let $c: V(P) \rightarrow \{0, 1, 2, 3\}$ be a nonsingular $4$-colouring of $P$ and suppose that $P(i, j)$ is a subgraph of $P$ which is induced on the vertices coloured $i$ and $j$ by $c$.
 Without loss of generality we may assume that~$v$ is coloured by $3$ and the edge $g_{q}$ is coloured by $(3, q)$. Let us denote

\begin{align}
 P^{0} : = P(3, 0) \cup P(1, 2), \  P^{1} : = P(3, 1) \cup P(0, 2)\  \hbox{ and } \ P^{2} : = P(3, 2) \cup P(0, 1).
\end{align}
Certainly, $\{E( P^{0}), E( P^{1}), E( P^{2})\}$ is a partition of the edge set of $P$. An edge (a subgraph) in $P$ is said to be of \textit{$q$-class} if this edge (any edge of this subgraph, respectively) belongs to $P^q$.  Certainly, the edge $g_{q}$ is of $q$-class, for $q \in Z_3$. Since $c$ is nonsingular the following proposition is satisfied.
\begin{proposition}\label{prop1.2}
If three edges in $P$ having a common vertex (say $v$) are successive edges w.r.t. counter-clockwise orientation around $v$, then they belong to successive classes ($P^{q}$, $P^{q+1}$, $P^{q+2}$, for some $q \in Z_3$). 
\end{proposition} 
From the above proposition, it follows that  $P^{q}$ is a spanning subgraph of $P$, for $q \in Z_3$. Moreover,  each edge of $q$-class belongs to only one maximal path of $q$-class (with ends of degree $3$ in $P$) or it belongs to only one cycle of $q$-class,  for $q \in Z_{3}$. Since $P$ has four vertices of degree~$3$, there are two maximal paths of $q$-class. Further, there exists a drawing of $P$ (denoted $P_{q}$, for $q \in Z_3$) satisfying the following conditions: The edges $g_{0}$,  $g_{1}$,  $g_{2}$ have counter-clockwise orientation around the common vertex in the triangulation $P_{q}$.  The triangulation $P_{q}$ consists of a maximal path of $q$-class (called \textit{inner path}) containing the edge $g_{q}$, which is situated on a line. This path has length $M(q)$ and is surrounded by $K(q) - 1$ cycles of $q$-class with the same length $2M(q)$. Finally, there is another maximal path of $q$-class with the length $M(q)$ (called \textit{exterior path}) around the outside of the last cycle (see Figs. 1 and 2).  

Since $P$ is $3$-connected the isomorphism between $P$ and $P_{q}$, for $q \in Z_3$, is combinatorial (Diestel \cite[p.93]{flo2}). Hence it is an orientation-preserving isomorphism.  
 Gr\"{u}nbaum and Motzkin  \cite{flo7} proved (in dual terms), without using the result of Fisk, that there exists the drawing $P_q$ of $P$, for $q \in Z_3$. By the definitions of $K(q)$ and $M(q)$ we obtain the following Gr\"{u}nbaum and Motzkin \cite [Lemma 2]{flo7} result  
\begin{proposition}\label{prop1.3} [Gr\"{u}nbaum and Motzkin] 
$$|P| = 2K(q)M(q)+2,\  \hbox{for}\  q \in Z_{3}.$$
\end{proposition}

Note that the exterior path  of $P_q$ may be situated at many different positions. Florek \cite[Definition 2.2]{flo5} defined an integer  $0\leqslant S^+ (q) < M(q)$  ($0 < S^- (q) \leqslant M(q)$) which determines the  position of this path (see also Definition \ref{def2.2} and Corollary \ref{cor2.1}  in chapter~$2$).  The vector $(K(q), M(q), S^+(q))$ is called the \textit{index-vector} of $P_{q}$ (or an index-vector of~$P$), for $q \in Z_3$, and the set $O(P) = \{(K(q), M(q), S^+(q): q \in Z_{3}\}$ is called the \textit{orbit} of $P$. 
\bigskip

Florek  \cite[Theorem 3.1 and Theorem 3.2]{flo5} introduced arithmetic equations which allow to calculate the index-vector $(K(q+1), M(q+1), S^+(q+1))$ of $P_{q+1}$ by the index-vector  $(K(q), M(q), S^+(q))$ of $P_{q}$, for $q \in Z_3$. The results are summarized in the following theorem:

\begin{theorem}\label{theorem1.1} Let $(K(q), M(q), S^+(q))$ and $(K(q+1), M(q+1), S^+(q+1))$, $q \in Z_{3}$, be index-vectors of $P \in \mathcal{P}$. Suppose that  $aM(q) - bS^{+}(q) = d > 0$, where $a, b$ are positive integers, $b \leqslant \afrac{M(q)}{d}$ and $d = gcd(S^+ (q), M(q))$. Then we have:
\bigskip

\begin{enumerate}
 \item[$(1)$]   $S^{+} (q+1)  \equiv  b K(q) - K (q+1) \pmod{ M(q+1)}$,
\item[$(2)$] $K(q+1) = gcd(S^+ (q), M(q))$,
\item[$(3)$] $K(q+ 1) M(q+ 1) = K(q) M(q)$.
\end{enumerate}
\end{theorem}

\noindent
The above theorem implies that drawings $P_{1}$, $P_{2}$, $P_{3}$ of the  triangulation $P \in \mathcal{P}$ cyclically determine one another. 
Theorem \ref{theorem1.1} yields also  the following proposition:
\begin{proposition}\label{proposition1.4} The following conditions are satisfied:
\begin{enumerate}
\item[$(1)$] any two triangulations of $\mathcal{P}$ are equivalent up to orientation-preserving isomorphism if and only if they have the same orbit,
\item[$(2)$] if two orbits of triangulations of $\mathcal{P}$ are different, then they are disjoint,
\item[$(3)$] an orbit of any triangulation of $\mathcal{P}$ is of order  $1$ or $3$.
\end{enumerate}
\end{proposition}
Suppose that $P$ has an orbit of degree $1$. It follows from Theorem\ref{theorem1.1} and Corollary~\ref{cor2.1} (see chapter $2$) that for every vertex $v$ of degree $3$ in $P$ there exists an automorphism~$\phi$ of the triangulation $P$ preserving orientation such $\phi(v) = v$ and $\phi(P^{q}) = P^{q+1}$, for $q \in Z_{3}$. Florek  \cite[Theorem 4.1]{flo5} characterized orbits of order $1$ of triangulations in $\mathcal{P}$ in the following way: 

\begin{proposition}\label{prop1.6} Let $k, m, s$ be integers, $k, m > 0$ and $0 \leqslant s < m$. Then,
$\{(k, m, s)\}$  is the orbit of some triangulation in $\mathcal{P}$ if and only if $(k, m, s) = (k, kz, kx)$, where integers $0 \leqslant x < z$ are solutions of the Diophantine equation $x^{2} + x + 1 = yz$. 
\end{proposition}
\bigskip

\noindent
Schinzel \cite[Appendix]{flo5} and \cite{flo11} has found formulas for all integers  $0 \leqslant x < z$ and $y$ which are solutions of the Diophantine equation $x^2+x+1 = yz$. 
\bigskip

\noindent
Let $t(n)$, $n \in N$, be the number of solutions of the congruence $x^{2} + x + 1  \equiv 0 \pmod{n}$ for $0 \leqslant x < n$. Mohar \cite[Proposition 3.2]{flo9} proved the following:
\begin{proposition}\label{prop1.7} [Mohar]
Let $n = 3^{\alpha}p^{\alpha_{1}}_{1}\ldots p^{\alpha_{j}}_{j}$ ($j \geqslant 0, \alpha \geqslant 0)$ be the decomposition of $n$ into primes. 
If  for $i = 1, 2, \ldots, j$, $p_{i} \equiv 1 \pmod{ 3}$, and $\alpha \in \{0, 1\}$, then $t(n) = 2^{j}$. In any other case $t(n) = 0$.
\end{proposition}
\bigskip

Let $\bar{P}$ be a mirror reflection of $P \in \mathcal{P}$. The union  $O(P) \cup O(\bar{P})$ is called a code of $P$. We prove the following proposition (in chapter $2$) 

\begin{proposition}\label{prop1.5} The following conditions are satisfied:
\begin{enumerate}
\item[$(1)$] any two triangulations of $\mathcal{P}$ are isomorphic  if and only if they have the same code,
\item[$(2)$] if two  codes of triangulations of $\mathcal{P}$ are different, then they are disjoint,
\item[$(3)$] a code of any triangulation of $\mathcal{P}$  is of order $j$, where $j = 1, 2, 3$ or $6$.
\end{enumerate}
\end{proposition}
\bigskip

There is a simple characterization of triangulations belonging to $\mathcal{P}$ that have codes of order $1$. Namely, by Lemma \ref{lemma5.2}, $P \in \mathcal{P}$ has a code of order $1$ iff  $O(P) = \{(\omega, \omega, 0)\}$ or  $O(P) = \{(\omega, 3\omega, \omega)\}$, for some positive integer $\omega$.
\bigskip

We say that $P$ is \textit{symmetric} if $P$ and $\bar{P}$ have the same  orbit. There is a simple characterization of symmetric triangulations belonging to $\mathcal{P}$. Let $\{(K(q), M(q), S^+(q): q \in Z_{3}\}$ be an orbit of $P \in \mathcal{P}$. It follows from Lemma \ref{lemma5.1}, that $P$ is symmetric if and only if $S^{+}(q) = M(q) - S^{-}(q)$, for some $q \in Z_{3}$. 
\bigskip

We can provide a characterization of triangulations belonging to $\mathcal{P}$ whose codes are of degree $2$ or $3$. Namely, by Lemma \ref{lemma6.0}, $P$ has a code of order $2$ if and only if $P$ has an orbit of order $1$ different than $\{(\omega, \omega, 0)\}$ and $ \{(\omega, 3\omega, \omega)\}$, for any integer~$\omega$. Moreover, by  Lemma \ref{lemma7.1}, $P$ has a code of order $3$ if and only if $P$ is symmetric and it does not have an orbit of the form  $\{(\omega, \omega, 0)\}$ and $ \{(\omega, 3\omega, \omega)\}$, for any positive integer~$\omega$.  
\bigskip

Let $\theta(n)$ be the number of all divisors of a positive integer $n$ and suppose that $\sigma(n)$ is the sum of all divisors of $n$. Let $\mathcal{P}_{n}$ denote the family of all triangulations of order $2n + 2$ in $\mathcal{P}$. We calculate  the number of symmetric triangulations of $\mathcal{P}_{n}$ (see Lemma~\ref{lemma5.3}). Using Propositions \ref{prop1.6}-\ref{prop1.7}, we enumerate  all triangulations of $\mathcal{P}_{n}$ which have codes of order~$2$ and of order $1$ (see Lemma \ref{lemma6.1}). Finally, we  prove the following theorem (in Chapter $6$): 
\begin{theorem}\label{theo1.1}
Let  $n = 2^{l}3^{\alpha}p^{\alpha_{1}}_{1}\ldots p^{\alpha_{u}}_{u}q^{\beta_{1}}_{1}\ldots q^{\beta_{w}}_{w}$ $(l, \alpha, \alpha_{i}, \beta_{i}\geqslant 0)$ be the decomposition of $n$ into primes such that $p_{i} \equiv 1 \pmod{ 3}$, for $i = 1,\ldots, u$, and $q_{i} \equiv 2 \pmod{ 3}$, for $i = 1,\ldots, w$. 
Let $d(n)$ be the number of (non-isomorphic) triangulations of $\mathcal{P}_{n}$.
\\[4pt]  
Then $d(1) = 1$ (tetrahedron is the only triangulation of $\mathcal{P}_{1}$).   \\[4pt] 
If $n > 1$ and $l = 0$, then,
\begin{align}
 d(n) = \frac{\sigma(n) +  3\theta(n) + 2 \theta^{\star}(n)}{6} -1.
\end{align}
If $n > 1$ and $l > 0$, then

  \begin{align}
  d(n) = \frac{\sigma(n) + 3(2l - 1)\theta(\afrac{n}{2^{l}}) + 2 \theta^{\star}(n)}{6} - 1 ,
\end{align}
\\[2pt]
where
\[
\begin{array}{ll}
\theta^{\star} (n) = 
\left\{
\begin{array}{ll}
\theta(p^{\alpha_{1}}_{1}\ldots p^{\alpha_{u}}_{u}), &\hbox{if} \ \l \ \hbox{and} \ \beta_{i} \ \hbox{are even}, \ \hbox{for every} \ i = 1, \ldots,  w, 
\\
 0, &\hbox{if}\ \ l \ \hbox{or} \ \beta_{i} \ \hbox{is odd}, \ \hbox{for some}\ i = 1, \ldots,  w.
\end{array}
\right.
\end{array}
\]
\end{theorem}

\begin{figure}
\centering

\begin{tikzpicture}[xscale=0.85, yscale=0.7361]

\coordinate  (A) at (0,0);
\coordinate  (B) at (6,0);

\coordinate  (E) at (1,0);
\coordinate  (F) at (2,0);
\coordinate  (G) at (3,0);
\coordinate  (H) at (4,0);
\coordinate  (J) at (5,0);

\coordinate (C) at (3.5,1);
\coordinate (K) at (3.5,2);
\coordinate (L) at (3.5,3);
\coordinate (M) at (3.5,4);

\coordinate (D) at (2.5,-1);
\coordinate (N) at (2.5,-2);
\coordinate (P) at (2.5,-3);

\draw (A)--node[near start,above,xshift=-7mm,yshift=-1mm]{$g_0$}(B);
\draw (C)--(M);
\draw (D)--(P);

\draw (A)--node[near start,above,xshift=-6mm,yshift=-5mm]{$g_1$}(M)--(E)--(L)--(F)--(K)--(G)--(C)--(H)--(K)--(J)--(L)--(B)--(M);
\draw (A)--node[near start,above,xshift=-4mm,yshift=-2mm]{$g_2$}(P)--(E)--(N)--(F)--(D)--(G)--(N)--(H)--(P)--(J);

\draw (J) .. controls (7,-2) and (9,3) .. (M);
\draw (M) .. controls (-2,3) and (-2,-2) .. (P);

\filldraw[black]
(A) circle (2pt)
(B) circle (2pt)
(C) circle (2pt)
(D) circle (2pt);

\draw (A) node[anchor=east]{$A$};
\draw (B) node[anchor=west]{$B$};
\draw (C) node[anchor=north, yshift=-1mm, xshift=-1pt]{$C$};
\draw (D) node[anchor=south, yshift=1mm, xshift=0.1mm]{$D$};

\end{tikzpicture}

\caption{The triangulation $P_0$ with the index-vector $(1,6,3)$, $S^{-}(0) = 4$.}
\end{figure}

\begin{figure}
\centering

\begin{tabular}[b]{c@{\hspace{-8mm}}c}

\begin{tikzpicture}[xscale=0.85, yscale=0.7361]

\coordinate (O) at (0,0);
\coordinate  (A) at (-1,0);
\coordinate  (B) at (1,0);

\coordinate (C) at (2.5,2.5);
\coordinate (D) at (-2.5,-2.5);
\coordinate (CC) at (3.5,3.5);

\coordinate  (E) at (-2,0);
\coordinate  (F) at (-2.5,0);
\coordinate  (G) at (2,0);
\coordinate  (H) at (2.5,0);

\coordinate (I) at (-0.5,1);
\coordinate (J) at (0.5,1);
\coordinate (K) at (-0.5,-1);
\coordinate (L) at (0.5,-1);
\coordinate (M) at (0,2);
\coordinate (N) at (0,-2);

\draw (A)--node[near start,below,xshift=-2mm,yshift=2.1mm]{$g_0$}(N)--(L)--(K)--(O)--node[above,xshift=-0.1mm,yshift=-1mm]{$g_1$}(A)--node[above,xshift=-2mm,yshift=-1.9mm]{$g_2$}(I)--(L)--(B)--(O)--(J)--(I)--(M)--(B);
\draw (I) .. controls (-2,1) and (-2,-1) ..(K);
\draw (J) .. controls (2,1) and (2,-1) ..(L);

\draw (I) .. controls (-1,1) and (-2,0.5) ..(F);
\draw (K) .. controls (-1,-1) and (-2,-0.5) ..(F);
\draw (J) .. controls (1,1) and (2,0.5) ..(H);
\draw (L) .. controls (1,-1) and (2,-0.5) ..(H);

\draw (H) .. controls (2.5,1) and (1,2) .. (M);
\draw (H) .. controls (2.5,-1) and (1,-2) .. (N);
\draw (F) .. controls (-2.5,1) and (-1,2) .. (M);
\draw (F) .. controls (-2.5,-1) and (-1,-2) .. (N);

\draw (CC) -- (H) -- (C) -- (CC) -- (M) -- (C);
\draw (N) -- (D) -- (F);

\draw (F) .. controls (-2.5,3) and (3,4) .. (CC);
\draw (N) .. controls (4,-2.5) and (4,3.5) .. (CC);
\draw (D) .. controls (4,-4.5) and (5,3.5) .. (CC);

\filldraw[black]
(A) circle (2pt)
(B) circle (2pt)
(C) circle (2pt)
(D) circle (2pt);

\draw (A) node[anchor=east]{$A$};
\draw (B) node[anchor=west]{$B$};
\draw (C) node[anchor=west]{$C$};
\draw (D) node[anchor=east]{$D$};

\end{tikzpicture}
 &

\begin{tikzpicture}[xscale=0.85,yscale=0.7361]

\coordinate (O) at (0,0);
\coordinate  (A) at (-1,0);
\coordinate  (B) at (0,-2);
\coordinate (C) at (2,0);
\coordinate (D) at (1,2);

\coordinate  (BB) at (0,-3);
\coordinate (DD) at (1,3);

\coordinate  (E) at (1,0);

\coordinate  (F) at (-0.5,1);
\coordinate  (G) at (0.5,1);
\coordinate  (H) at (1.5,1);

\coordinate  (I) at (-0.5,-1);
\coordinate  (J) at (0.5,-1);
\coordinate  (K) at (1.5,-1);

\coordinate (fake1) at (0,3.99);
\coordinate (fake2) at (0,-4.5);

\draw (A)--node[near start,above,xshift=-1mm,yshift=-0.6mm]{$g_0$}(F)--(O)--node[above,xshift=-0.1mm,yshift=-1mm]{$g_2$}(A)--node[below,xshift=-1.5mm,yshift=1.9mm]{$g_1$}(I)--(O)--(G)--(E)--(O)--(J)--(E)--(H)--(C)--(E)--(K)--(C);
\draw (DD)--(F)--(H)--(DD)--(G)--(D)--(H)--(DD)--(D);
\draw (I)--(BB)--(K)--(I)--(B)--(J)--(BB)--(B);

\draw (F) .. controls (-2,1) and (-2,-1) ..(I);
\draw (H) .. controls (3,1) and (3,-1) ..(K);

\draw (F) .. controls (-3,1) and (-2.5,-3) ..(BB);
\draw (K) .. controls (4,-1) and (3,3) ..(DD);
\draw (BB) .. controls (4,-2) and (4,3) ..(DD);

\filldraw[black]
(A) circle (2pt)
(B) circle (2pt)
(C) circle (2pt)
(D) circle (2pt);

\draw (A) node[anchor=east]{$A$};
\draw (B) node[anchor=south,yshift=1mm]{$B$};
\draw (C) node[anchor=west]{$C$};
\draw (D) node[anchor=north,yshift=-1mm]{$D$};

\filldraw[white] (fake1) circle (0.1pt);
\filldraw[white] (fake2) circle (0.1pt);

\end{tikzpicture}
\end{tabular}

\caption{The triangulations $P_1$ with the index-vector $(3,2,0)$, $S^{-}(1) = 1$, and $P_2$ with the index-vector $(2,3,1)$, $S^{-}(2) = 3$.}
\end{figure}

In the last section we will present application of equation (2) from Theorem \ref{theo1.1}. 

A positive integer $n$ is said to be \textit{perfect} if and only if $\sigma(n) = 2n$. Euler \cite[p. 19]{flo12}) established a multiplicative structure of odd perfect numbers:
\begin{theorem}\label{theorem1.3} $[Euler]$
If $n$ is an odd perfect number, then
\begin{align*} 
n = \pi^{4s+1}M^{2},\hbox { where }
\pi \hbox{ is prime }, gcd(\pi,M) = 1 \hbox{  and } \pi  \equiv  n \equiv 1 \pmod{4}.
\end{align*}
\end{theorem}
Using  equality (2) of Theorem \ref{theo1.1}, we obtain the following strengthening of the Euler theorem: 
\begin{theorem} \label{theo1.4} If $n$ is an odd perfect number and $3 \nmid n$, then
\begin{align*} 
 n = \pi^{12s +1}M^{2} \hbox{ or }  n = \pi^{12s +9}M^{2},& \hbox { where }
\pi \hbox{ is prime}, \ gcd(\pi,M) = 1 
\cr
&\hbox{ and } \pi  \equiv  n \equiv 1 \pmod{12}. 
 \end{align*}
\end{theorem}
\section{Basic definitions}
\bigskip

Let $P \in \mathcal{P}$ and suppose that $P^{0}$, $P^{1}$, $P^{2}$ are subgraphs of $P$ defined by condition~$(1)$ in Introduction.  Hence, they satisfy Proposition \ref{prop1.2}. Recall that $P^{q}$ contains two maximal paths of $q$-class both with the same length $M(q)$, and $K(q) - 1$ cycles of $q$-class both with the same length $2M(q)$, for every $q \in Z_3$. 

The following Definition \ref{def2.1}, Definition \ref{def2.2}, Lemma \ref{lem2.1} and Theorem \ref{theo2.1} were given and proved,  by Florek  \cite{flo5}. 
\begin{definition}\label{def2.1} 
Let $[A, q]= v_0 v_1 \ldots \ v_{M(q)}$ be a maximal path of $q$-class in $P$. We may assume that it is a directed path such that $A = v_0$ is its initial and $A_q = v_{M(q)}$ is its terminal vertex, both of degree $3$ in $P$. An edge $e$ adjacent to the directed path $[A, q]$ is called a~\textit{left branch} of the path if it is branching off from $[A, q]$ to the left. More precisely, $(v_{i} v_{i+1}, e)$ for $0\leqslant i < M(q)$ ($(e, v_{i-1}v_i)$ for $0 < i \leqslant M(q)$) is a pair of counter-clockwise successive edges incident with the vertex $v_i$. Otherwise, it is called  a~\textit{right branch} of the path. We put
$$[ A, q] (e) = \left\{\begin{array}{ll} i & \quad \mbox{if $e$ is a~left branch of\/ $[A, q ]$ incident with $v_i$},
\\[4pt]
 2M(q) - i & \quad \mbox{if $e$ is a~right branch  of\/ $[A, q ]$ incident with $v_i$}.
\end{array}\right.
$$
\end{definition}
\begin{remark}\label{rem2.1}
Notice that $e$ is a left branch of the path $[A,q]$ if and only if it is a right branch of the path $[A_q, q]$. Moreover, $\left| [A_q, q ] (e) - [ A, q ] (e) \right| = M(q)$.
\end{remark}
\begin{lemma}\label{lem2.1}
Let $A, C$ be ends of two different maximal paths of $q$-class. Suppose that $e,\hat{e}$ and $f,\hat{f}$ are pairs of end-edges  of  two minimal paths of $(q+1)$-class so that $e$ is incident with $A$, $\hat{e}$ is adjacent to the path $[C, q]$, $f$ is adjacent to the path $[A, q]$ and $\hat{f}$ is incident with $C$. Then 
$$
[A, q ] (f) = [C, q](\hat{e}).
$$
\end{lemma}
\begin{definition}\label{def2.2}
Let $A$, $C$ be ends of two different maximal paths of $q$-class in the triangulation $P$, and suppose $f$
\1{(}or $g$\1{)} is the first edge of the directed path $[C, q+1 ]$ \1{(}or $[C, q-1]$, respectively\1{)}
which is adjacent to the directed path $[A, q] = v_0 v_1 \ldots v_{M(q)}$.
$$S^+ (q) := \left\{\begin{array}{ll}
i,  & \quad \mbox{if $f$  is a~left branch of $[A, q ]$ incident with $v_i$},
\\[4pt]
 M(q) - i, & \quad \mbox{if $f$ is a~right branch of $[A, q ]$ incident with $v_i$}.
\end{array}\right.
$$
$$S^- (q) := \left\{\begin{array}{ll}
i, & \quad \mbox{if $g$  is a~left branch of $[A, q ]$ incident with $v_i$},
\\[4pt]
M(q) - i, & \quad \mbox{if $g$ is a~right branch of $[A, q ]$ incident with $v_i$}.
\end{array}\right.
$$
\end{definition}

\begin{corollary}\label{cor2.1} By Remark~\ref{rem2.1} and by Lemma \ref{lem2.1} definition of $S^+ (q)$ and $S^- (q)$ do not depend on the choice of ends of two different maximal paths of $q$-class. Hence definition of the index-vector of $P_q$, definition of the orbit and of the code of a triangulation $P \in\mathcal{P}$ do not depend on the choice of a vertex of degree $3$ in $P$.
\end{corollary}

Florek \cite[Theorem 2.1]{flo5} shows that $ S^+(q)$ and $ S^-(q)$ determine one another:
\begin{theorem}\label{theo2.1} Let $(K(q), M(q), S^{+}(q))$ be an index-vector of $P$, for some $q \in Z_3$.
$$
S^- (q) - S^+(q) \equiv K(q) \pmod {M(q)} .
$$
\end{theorem}

\begin{corollary}\label{cor2.2}
By definitions of $S^{+}(q)$ and $S^{-}(q)$, if $(K(q), M(q), S^+(q))$  is an index-vector of $P_q$, then $(K(q), M(q), M(q) - S^-(q))$ is an index-vector of $\bar{P}_q$. Hence, by Theorem \ref{theo2.1}, we have
\begin{enumerate}
\item[$(i)$] $P$ and $R$ have the same orbit if and only if $\bar{P}$ and $\bar{R}$ have the same orbit,
\item[$(ii)$] $P$ and $\bar{P}$ have orbits of the same order.
\end{enumerate}
\end{corollary}

\begin{proof} [\bf Proof of Proposition \ref{prop1.5}]
\bigskip

Let $P, R \in \mathcal{P}$ and suppose that $O(P), O(R), O(\bar{P}), O(\bar{R})$ denote orbits of $P, R, \bar{P}, \bar{R}$, respectively. Then, 
$O(P)\cup O(\bar{P})$, $O(R)\cup O(\bar{R})$ are codes of $P$ and $R$, respectively.
\bigskip

Case (a). $P$ and $R$ are isomorphic if and only if $P$ and $R$,  or $P$ and $\bar{R}$, are equivalent up to orientation-preserving isomorphism (because $P$ is $3$-connected). Hence, equivalently,  $O(P) = O(R)$,  or $O(P) =  O(\bar{R})$, by condition (1) of Proposition \ref{proposition1.4}. Thus, equivalently,  $O(P)\cup O(\bar{P}) = O(R)\cup O(\bar{R})$, by condition (i) of Corollary \ref{cor2.2}. 
\bigskip

Case (b). If codes of $P$ and $R$ are different, then, by condition (i) of Corollary \ref{cor2.2}, $O(P) \neq  O(R)$, $O(P) \neq O(\bar{R})$, $O(\bar{P})\neq O(R)$ and $O(\bar{P}) \neq O(\bar{R})$. Therefore,  by condition~(2) of Proposition \ref{proposition1.4}, $O(P)\cup O(\bar{P})$ and $O(R)\cup O(\bar{R})$ are disjoint.
\bigskip

Case (c). By condition (3) of Proposition \ref{proposition1.4}, $P$ has orbit of order $j$,  for $j = 1$ or~$3$. It follows from condition (ii) of Corollary \ref{cor2.2}, that $\bar{P}$ has orbit of order $j$. Hence, by condition (2) of Proposition \ref{proposition1.4}, $P$ has code of order $j$ or $2j$. Therefore, if $P$ is symmetric (is not symmetric) $P$ has code of order $1$ or $3$ (of order $2$ or $6$, respectively).

\end{proof}
\bigskip

\section{Non-simple triangulations}\bigskip
 
Let $\mathcal{P}^{*}_{n}$ be the family of all triangulations of order $2n + 2$ with all vertices of degree $3$ or $6$. 
\begin{proposition}\label{prop4.1} The following conditions are satisfied:
\begin{enumerate}
\item[$(1)$] each non-simple triangulation of $\mathcal{P}^{*}_{n}$,   has two non adjacent edges each of which has end-vertices of degree $3$, 
\item[$(2)$] if $G$ (or $G_1$) is a non-simple triangulation of $\mathcal{P}^{*}_{n}$ having an edge $cd$ ($c_{1}d_{1}$, respectively) with end-vertices of degree $3$, then there exists isomorphism $\sigma : G \rightarrow G_{1}$ such that $\sigma(c) = c_{1}$ and $\sigma(d) = d_{1}$.
\end{enumerate}
\end{proposition}
\begin{proof} Let $n > 1$ and suppose that $G \in\mathcal{P}^{*}_{n}$ is not simple. Then, $G$ contains a cycle $aba$ of length $2$. The cycle $aba$ determines bounded (say $U'$) and unbounded (say $U''$) regions on the plane. Since $P$ is a triangulation, both $U$ and $U''$ contain at least one vertex. Notice that if  $a$ is adjacent with only one vertex (say $e$) belonging to $U'$, then there are two triangles with the same set of vertices $\{a, b, e\}$ which is a contradiction because, by Proposition \ref{prop1.1}, $G$ has a nonsingular $4$-colouring. Hence, $a$ and~$b$ have two common neighbours belonging to $U'$ and two common neighbours belonging to $U''$. If we delete vertices belonging to $U'$ (or $U''$) and one of the edges $ab$ we obtain a triangulation~$G'$ (or~$G''$, respectively) belonging to $\mathcal{P}^{*}_{m}$, where $1 < m <n$, or $G' = K_{4}$ ($G''= K_{4}$, respectively).  If, by induction, $G'$ (or $G''$) has a pair of non-adjacent edges with end-vertices of degree $3$, then one of them is different from $ab$. Hence, $G$ has also a pair of non-adjacent edges with end-vertices of degree $3$. Therefore, condition $(1)$ holds.

Let $G$ and $G_{1} \in\mathcal{P}^{*}_{n}$ be not simple. Then, by condition $(1)$, $G$ (or $G_1$) has an edge $cd$ (or $c_{1}d_{1}$) with end-vertices of degree $3$. Hence, there are two triangles $acd$ and $bcd$ ($a_{1}c_{1}d_{1}$ and $b_{1}c_{1}d_{1}$, respectively) such that $aba$ ($a_{1}b_{1}a_{1}$, respectively) is a cycle of length~$2$. If we delete vertices $c, d$ ($c_{1}, d_{1}$) and one edge with end-vertices $a, b$ ($a_{1}, b_{1}$) from~$G$ ($G_{1}$, respectively) we obtain a triangulation $G'$ (or $G'_{1}$, respectively) belonging to~$\mathcal{P}^{*}_{n-1}$, for $n > 2$, or $G' = G'_{1} = K_{4}$, for $n = 2$. By induction, there exists an isomorphism $\sigma' : G' \rightarrow G'_{1}$ such that $\sigma'(a) = a_{1}$ and $\sigma'(b) = b_{1}$. Certainly, we may extend it to the isomorphism $\sigma : G \rightarrow G_{1}$ such that $\sigma (c) = c_{1}$ and $\sigma(d) = d_{1}$. Hence, condition $(2)$ holds.
\end{proof}
\begin{remark}\label{rem4.1}
Let $X_{n} = \{(k, m , s)\in \mathbb{Z}^{3}: \ 1 \leqslant m \leqslant n, \  0 \leqslant s  < m \ \hbox{and} \ km = n\}$. 
A vector $(k, m , s) \in X_{n}$ is called \textit{proper} if it is different from $(n,1,0)$, $(1, n, n-1)$ and $(1,n, 0)$, for $n > 1$. 
Notice that, by condition $(1)$ of Proposition \ref{prop4.1}, each vector of $X_{n}$ is proper if and only if it is an index-vector of some triangulation in $\mathcal{P}_{n}$ (of order $2n+2$).  We may say that $\{ (n,1,0), (1, n, n-1), (1,n, 0)\}$ is the code of the non-simple graph in~$\mathcal{P}^{*}_{n}$, for $n > 1$.
\end{remark}

\section{Symmetric triangulations}
\bigskip

Let $P \in \mathcal{P}$ and suppose that  $P_{q}$, for $q \in Z_{3}$, is a drawing of~$P$ with an index-vector $(K(q), M(q), S^{+}(q))$. We recall that $P_{q}$ has two maximal paths of $q$-class, called the inner and the exterior path of $P_{q}$. The inner path is situated on a line  (say $l_{q}$). Let $\bar{P}_{q}$ be obtained from $P_{q}$ after a transformation of symmetry with respect to $l_{q}$. Let $v_{0}v_{1}\ldots v_{M(q)}$ be the exterior path of $P_q$. If  $P$ is not tetrahedron, then $M(q) > 1$ because~$P$ is simple. 

 We say that $P_q$ is a \textsl{mirror symmetric} drawing of $P$ if one of the following (equivalent) conditions is satisfied: 
\begin{enumerate}
\item[$(i)$] $P_q$ and $\bar{P}_{q}$ have the same index-vector,
\item[$(ii)$] $S^{+}(q) = M(q) - S^{-}(q)$,
\item[$(iii)$] $l_{q}$ is the axis of symmetry of $P_{q} - v_{1}v_{2}\ldots v_{M(q)-1}$.
\end{enumerate}

Let $P_{q}$  be a mirror symmetric drawing of some triangulation $P \in \mathcal{P}$. Then, $P_{q}$ is called of \textsl{ type} $1$  (of \textsl{ type} $2$) if ends of the exterior path of $P_{q}$ are situated on a line perpendicular to the line $l_q$  (on the line $l_q$, respectively).
\bigskip

\begin{example} \label{example5.1} Triangulations of Figs 1 and 2  are not mirror symmetric drawings of any triangulation of $\mathcal{P}$.  Triangulations of Figs 3 and 4 are mirror symmetric drawings of some triangulations in $\mathcal{P}$. Triangulation of Fig. 3 on the left (on the right) is of type $1$ (of type~$2$, respectively).  Similarly, triangulation of Fig. 4 on the left (on the right) is of type $1$ (of type $2$, respectively).
\end{example}
It is not difficult to check the following remark. 
\begin{remark}\label{rem5.1}  Let $P_{q}$  be a mirror symmetric drawing of some triangulation $P \in \mathcal{P}$ with an index-vector $(K(q), M(q), S^{+}(q))$.
\begin{enumerate}
\item[$(1)$] if $P_{q}$  is of type $1$, then $K(q)$ and $M(q)$ are of the same parity  (see Figs 3, 4),
\item[$(2)$] if $P_{q}$  is of type $2$, then $K(q)$ is even. 
\end{enumerate}
\end{remark}
We recall that $\mathcal{P}_{n}$ denote the family of all triangulations of order $2n + 2$ in $\mathcal{P}$. It is easy to check the following remark. 
\begin{remark}\label{rem5.2} Let $n > 1$. 
\begin{enumerate}
\item[$(1)$] If  $(k, m)$ is a pair of positive integers of the same parity, $m > 1$  and $n = km$, then there exists only one mirror symmetric drawing of some triangulation in $\mathcal{P}_{n}$, which is of type $1$ and $(k, m, s)$ is its index vector, for some integer $s$.
\item[$(2)$] If $(k, m)$ is a pair of positive integers, $k$ is even, $m > 1$ and $n = km$, then there exists only one mirror symmetric drawing of some triangulation belonging to $\mathcal{P}_{n}$, which is of type $2$ and $(k, m, s)$ is its index vector, for some integer $s$.
\end{enumerate}
\end{remark}
\begin{lemma}\label{lemma5.1} Let $P \in \mathcal{P}$. The following conditions are equivalent:
\begin{enumerate} 
\item[$(1)$] $P$ is symmetric,
\item[$(2)$] $P_q$  is a mirror symmetric drawing of $P$, for some $q \in Z_{3}$.
\end{enumerate}
\end{lemma} 
\begin{proof} $(1) \Rightarrow (2)$. If $P$ is symmetric, then $P$ and $\bar{P}$ have the same orbit. Thus, we have
\begin{enumerate}
\item[$(3)$]
$\{(K(q), M(q), S^{+}(q)): q \in Z_{3}\} = \{(K(q), M(q), M(q) - S^{-}(q)): q \in Z_{3}\}$.
\end{enumerate}
Assume that $(K(1), M(1), S^{+}(1)) = (K(2), M(2), M(2) - S^{-}(2))$.  Then, by Theorem \ref{theo2.1}
$$S^{+}(1) +  S^{-}(2) \equiv S^{-}(1) - K(1) + S^{+}(2) + K(2) = S^{-}(1) +  S^{+}(2)\pmod {M(1) = M(2)}.$$ 
Hence, $S^{-}(1) +  S^{+}(2) = M(1)$ because $S^{+}(1) + S^{-}(2) = M(2)$.  So $(K(2), M(2), S^{+}(2)) = (K(1), M(1), M(1) - S^{-}(1))$. Thus, by equality $(3)$, 
$$(K(3), M(3), S^{+}(3)) = (K(3), M(3), M(3) - S^{-}(3)).$$ Hence, $P_{3}$ and $\bar{P_{3}}$ have the same index-vector and condition $(2)$ holds. 

$(2) \Rightarrow (1)$. If condition (2) is satisfied, then, $P_q$ and $\bar{P}_{q}$ have the same index-vector, for some $q \in Z_{3}$. Hence, by condition $(2)$ of Proposition \ref{proposition1.4}, $P$ and  $\bar{P}$ have the same orbit and condition (1) holds. 
\end{proof}

\begin{lemma}\label{lemma5.2} Let $P \in \mathcal{P}$. The following conditions are equivalent:
\begin{enumerate}
\item[$(1)$] $P_{q_1}$ and $P_{q_2}$ are mirror symmetric drawings of $P$, for some $q_{1}\neq q_{2} \in Z_{3}$,
\item[$(2)$] $P$ has an orbit of the form $\{(\omega, \omega, 0)\}$ or  $\{(\omega, 3\omega, \omega)\}$, for some positive integer~$\omega$,
\item[$(3)$] $P$ has a code of order $1$.
\end{enumerate}
\end{lemma} 
\begin{proof} $(1) \Rightarrow (2)$.  Florek \cite[Theorem $4.2$]{flo5} proved that conditions (1) and (2) are equivalent.

$(2) \Rightarrow (3)$. If $\{(\omega, \omega, 0)\}$ (or $\{(\omega, 3\omega, \omega)\}$) is an orbit of a triangulation $P$, for some positive integer $\omega$, then, by Theorem \ref{theo2.1}, $\{(\omega, \omega, 0)\}$ ($\{(\omega, 3\omega, \omega)\}$, respectively) is the orbit of the triangulation $\bar{P}$. Hence, $P$ is symmetric and it has a code of order $1$.

$(3) \Rightarrow (1)$. If $P$ has a code of order $1$, then $P$ is symmetric and it has an orbit of order~$1$. Hence, $P_q$ and $\bar{P}_{q}$ have the same index-vector, for every $q \in Z_3$. Thus $P_q$ is a mirror symmetric drawing of $P$, for every $q \in Z_3$.
\end{proof}

\begin{figure}

\begin{tikzpicture}[xscale=0.5, yscale=0.4]

\coordinate  (a) at (-6.5,0);

\coordinate  (b) at (-3,0);
\coordinate  (c) at (-1,0);
\coordinate  (d) at (1,0);
\coordinate  (e) at (3,0);

\coordinate  (f) at (6.5,0);

\coordinate  (b3) at (-2,2);
\coordinate (c2) at (0,2);
\coordinate (e3) at (2,2);

\coordinate (b4) at (-1,4);
\coordinate (e4) at (1,4);

\coordinate (g) at (0,6);
\coordinate (g1) at (0,4.7);

\coordinate  (b1) at (-2,-2);
\coordinate (c1) at (0,-2);
\coordinate (e1) at (2,-2);

\coordinate (b2) at (-1,-4);
\coordinate (e2) at (1,-4);

\coordinate (h) at (0,-6);
\coordinate (h1) at (0,-4.7);

\draw (b)--(e)--(h)--(b)--(g)--(e);

\draw (b3)--(c)--(e4)--(b4)--(d)--(e3)--(b3);

\draw (b1)--(c)--(e2)--(b2)--(d)--(e1)--(b1);

\draw (b4)--(g1)--(g);
\draw (e4)--(g1);

\draw (b2)--(h1)--(h);
\draw (e2)--(h1);

\draw (b1) .. controls (-5,-1) and (-5,1) .. (b3);

\draw (e1) .. controls (5,-1) and (5,1) .. (e3);

\draw (a) .. controls (-5,1.5) and (-4,1.7)  .. (b3);

\draw (a) .. controls (-5,2) and (-3,3.5)  .. (b4);

\draw (a) .. controls (-5,-1.5) and (-4,-1.7)  .. (b1);

\draw (a) .. controls (-5,-2) and (-3,-3.5)  .. (b2);

\draw (f) .. controls (5,1.5) and (4,1.7)  .. (e3);

\draw (f) .. controls (5,2) and (3,3.5)  .. (e4);

\draw (f) .. controls (5,-1.5) and (4,-1.7)  .. (e1);

\draw (f) .. controls (5,-2) and (3,-3.5)  .. (e2);

\filldraw[black]
(a) circle (2pt)
(b) circle (4pt)
(c) circle (2pt)
(d) circle (2pt)
(e) circle (4pt)
(f) circle (2pt)
(b1) circle (2pt)
(b2) circle (2pt)
(b3) circle (2pt)
(b4) circle (2pt)
(e1) circle (2pt)
(e2) circle (2pt)
(e3) circle (2pt)
(e4) circle (2pt)
(c1) circle (2pt)
(c2) circle (2pt)
(g) circle (2pt)
(h) circle (2pt)
(g1) circle (4pt)
(h1) circle (4pt);

\begin{scope}[shift = {(15, 0)}]

\coordinate  (a) at (-6.5,0);

\coordinate  (b) at (-3,0);
\coordinate  (c) at (-1,0);
\coordinate  (d) at (1,0);
\coordinate  (e) at (3,0);

\coordinate  (f) at (6.5,0);

\coordinate  (b3) at (-2,2);
\coordinate (c2) at (0,2);
\coordinate (e3) at (2,2);

\coordinate (b4) at (-1,4);
\coordinate (e4) at (1,4);

\coordinate (g) at (0,6);
\coordinate (g1) at (0,4.7);

\coordinate  (b1) at (-2,-2);
\coordinate (c1) at (0,-2);
\coordinate (e1) at (2,-2);

\coordinate (b2) at (-1,-4);
\coordinate (e2) at (1,-4);

\coordinate (h) at (0,-6);
\coordinate (h1) at (0,-4.7);

\coordinate (x) at (-5,0);
\coordinate (y) at (5,0);

\draw (b)--(e)--(e1)--(b1)--(b)--(b3)--(e3)--(e);

\draw (b3)--(c1)--(e3);
\draw (b1)--(c2)--(e1);

\draw (b1) .. controls (-4.5,-1) and (-4.5,1) .. (b3);

\draw (e1) .. controls (4.5,-1) and (4.5,1) .. (e3);

\draw (a) .. controls (-5,1.5) and (-4,1.7)  .. (b3);

\draw (a) .. controls (-5,-1.5) and (-4,-1.7)  .. (b1);

\draw (f) .. controls (5,1.5) and (4,1.7)  .. (e3);

\draw (f) .. controls (5,-1.5) and (4,-1.7)  .. (e1);

\filldraw[black]
(a) circle (2pt)
(b) circle (4pt)
(c) circle (2pt)
(d) circle (2pt)
(e) circle (4pt)
(f) circle (2pt)
(b1) circle (2pt)

(b3) circle (2pt)

(e1) circle (2pt)

(e3) circle (2pt)

(c1) circle (2pt)
(c2) circle (2pt);

\filldraw [black] (x) circle (4pt);

\draw (x) .. controls (-4.5,-1.3) and (-3.5,-1.5)  .. (b1);
\draw (x) .. controls (-4.5,1.3) and (-3.5,1.5)  .. (b3);
\draw (x)--(a);

\filldraw [black] (y) circle (4pt);

\draw (y) .. controls (4.5,-1.3) and (3.5,-1.5)  .. (e1);
\draw (y) .. controls (4.5,1.3) and (3.5,1.5)  .. (e3);
\draw (y)--(f);

\end{scope}

\end{tikzpicture}

\caption{A triangulation $P_{q}$ with the index  $(K(q), M(q), S^{+}(q)) = (3, 3, 0)$ ($(2, 3, 2)$, respectively). $S^{+}(q) + S^{-}(q) = M(q)$.  The line containing the inner path of $P_q$ is the axis of symmetry of  $P_{q} - v_{1}v_{2}$, where  $v_{0}v_{1}v_{2}v_{3}$ is the exterior path of $P_q$.}
\end{figure}

\begin{figure}

\begin{tikzpicture}[xscale=0.5, yscale=0.4]


\coordinate  (a) at (1,4);
\coordinate  (aa) at (1,2.7);

\coordinate  (b) at (-3,0);
\coordinate  (c) at (-1,0);
\coordinate  (d) at (1,0);
\coordinate  (e) at (3,0);
\coordinate  (ee) at (5,0);

\coordinate  (f) at (1,-4);
\coordinate  (ff) at (1,-2.7);

\coordinate  (b3) at (-2,2);
\coordinate (c2) at (0,2);
\coordinate (e3) at (2,2);
\coordinate (ee3) at (4,2);

\coordinate (b4) at (-1,4);
\coordinate (e4) at (1,4);

\coordinate (g) at (0,6);
\coordinate (g1) at (0,4.7);

\coordinate  (b1) at (-2,-2);
\coordinate (c1) at (0,-2);
\coordinate (e1) at (2,-2);
\coordinate (ee1) at (4,-2);

\coordinate (b2) at (-1,-4);
\coordinate (e2) at (1,-4);

\coordinate (h) at (0,-6);
\coordinate (h1) at (0,-4.7);

\coordinate (x) at (-5,0);
\coordinate (y) at (7,0);

\draw (b)--(e)--(e1)--(b1)--(b)--(b3)--(e3)--(e);

\draw (b3)--(c1)--(e3);
\draw (b1)--(c2)--(e1);

\draw (e1)--(ee1)--(ee)--(ee3)--(e3);
\draw (ee1)--(e)--(ee3);
\draw (ee)--(e);

\draw (b1) .. controls (-4.5,-1) and (-4.5,1) .. (b3);

\draw (ee1) .. controls (6.5,-1) and (6.5,1) .. (ee3);

\filldraw[black]
(a) circle (2pt)
(b) circle (4pt)
(c) circle (2pt)
(d) circle (2pt)
(e) circle (2pt)
(ee) circle (4pt)
(f) circle (2pt)
(b1) circle (2pt)

(b3) circle (2pt)

(e1) circle (2pt)
(ee1) circle (2pt)
(e3) circle (2pt)
(ee3) circle (2pt)
(c1) circle (2pt)
(c2) circle (2pt)
(aa) circle (4pt)
(ff) circle (4pt);

\draw (c2)--(a)--(e3)--(aa)--(c2);
\draw (c1)--(f)--(e1)--(ff)--(c1);
\draw (a)--(aa);
\draw (f)--(ff);

\begin{scope}[shift = {(15, 0)}]
\coordinate  (a) at (-6.5,0);

\coordinate  (b) at (-3,0);
\coordinate  (c) at (-1,0);
\coordinate  (d) at (1,0);
\coordinate  (e) at (3,0);
\coordinate  (ee) at (5,0);

\coordinate  (f) at (8.5,0);

\coordinate  (b3) at (-2,2);
\coordinate (c2) at (0,2);
\coordinate (e3) at (2,2);
\coordinate (ee3) at (4,2);

\coordinate (b4) at (-1,4);
\coordinate (e4) at (1,4);

\coordinate (g) at (0,6);
\coordinate (g1) at (0,4.7);

\coordinate  (b1) at (-2,-2);
\coordinate (c1) at (0,-2);
\coordinate (e1) at (2,-2);
\coordinate (ee1) at (4,-2);

\coordinate (b2) at (-1,-4);
\coordinate (e2) at (1,-4);

\coordinate (h) at (0,-6);
\coordinate (h1) at (0,-4.7);

\coordinate (x) at (-5,0);
\coordinate (y) at (7,0);

\draw (b)--(e)--(e1)--(b1)--(b)--(b3)--(e3)--(e);

\draw (b3)--(c1)--(e3);
\draw (b1)--(c2)--(e1);

\draw (e1)--(ee1)--(ee)--(ee3)--(e3);
\draw (ee1)--(e)--(ee3);
\draw (ee)--(e);

\draw (b1) .. controls (-4.5,-1) and (-4.5,1) .. (b3);

\draw (ee1) .. controls (6.5,-1) and (6.5,1) .. (ee3);

\draw (a) .. controls (-5,1.5) and (-4,1.7)  .. (b3);

\draw (a) .. controls (-5,-1.5) and (-4,-1.7)  .. (b1);

\draw (f) .. controls (7,1.5) and (6,1.7)  .. (ee3);

\draw (f) .. controls (7,-1.5) and (6,-1.7)  .. (ee1);

\filldraw[black]
(a) circle (2pt)
(b) circle (4pt)
(c) circle (2pt)
(d) circle (2pt)
(e) circle (2pt)
(ee) circle (4pt)
(f) circle (2pt)
(b1) circle (2pt)

(b3) circle (2pt)

(e1) circle (2pt)
(ee1) circle (2pt)
(e3) circle (2pt)
(ee3) circle (2pt)
(c1) circle (2pt)
(c2) circle (2pt);

\filldraw [black] (x) circle (4pt);

\draw (x) .. controls (-4.5,-1.3) and (-3.5,-1.5)  .. (b1);
\draw (x) .. controls (-4.5,1.3) and (-3.5,1.5)  .. (b3);
\draw (x)--(a);

\filldraw [black] (y) circle (4pt);

\draw (y) .. controls (6.5,-1.3) and (5.5,-1.5)  .. (ee1);
\draw (y) .. controls (6.5,1.3) and (5.5,1.5)  .. (ee3);
\draw (y)--(f);

\end{scope}

\end{tikzpicture}

\caption{A triangulation $P_q$ with the index $(K(q), M(q), S^{+}(q)) = (2, 4, 1)$ ($(2, 4, 3)$, respectively). $S^{+}(q) + S^{-}(q) = M(q)$.   The line containing the inner path of $P_q$ is the axis of symmetry of  $P_{q}- v_{1}v_{2}v_{3}$, where  $v_{0}v_{1}v_{2}v_{3}v_{4}$ is the exterior path of $P_q$.} 
\end{figure} 

We recall that $\theta(n)$ is the number of all divisors of a positive integer $n$.
 \begin{lemma}\label{lemma5.3} 
 If $n$ is an odd integer, $n > 1$, then there exist $\theta(n) - 1$ symmetric triangulations (non-isomorphic) in $\mathcal{P}_n$. If  $n = 2^{l}(\afrac {n}{2^{l}})$, where~$l$ is a positive integer and $\afrac {n}{2^{l}}$ is odd, then there exist $(2l - 1) \theta (\afrac{n}{2^{l}}) - 1$ symmetric triangulations  in $\mathcal{P}_n$.
\end{lemma}
\begin{proof} 
Let  $\mathcal{S}_{n}$ be the set of all symmetric triangulations in $\mathcal{P}_n$. By Lemma \ref{lemma5.1}, $P \in \mathcal{S}_{n}$ if and only if there exists a drawing of $P$ which is a mirror symmetric drawing of $P$. Let $\mathcal{M}_{n}$ be the set of all index-vectors of mirror symmetric drawings of triangulations in $\mathcal{P}_{n}$. Hence, we may define a function $\omega: \mathcal{S}_{n}  \rightarrow \mathcal{M}_{n}: \omega(P)$ is an index-vector of a mirror symmetric drawing of $P \in \mathcal{S}_{n}$. If $P_{q_1}$ and $P_{q_2}$ are mirror symmetric drawings of~$P$, for $q_1 \neq q_2$, then, by  Lemma \ref{lemma5.2}, $P$ has code of order $1$. Hence, $P_{q_1}$ and $P_{q_2}$ have the same index-vector and the function $\omega$ is well defined.  Moreover, if $P_{q_{1}}$ and $R_{q_2}$ are mirror symmetric drawings of $P, R \in \mathcal{S}_{n}$ which have the same index-vector, then, by conditions (2) of Proposition \ref{proposition1.4}, $P$ and $R$ have the same orbit. Thus, by condition (1) of Proposition \ref{proposition1.4}, $P$ and $R$ are isomorphic. Hence,  the function $\omega: \mathcal{S}_{n}  \rightarrow \mathcal{M}_{n}$ is a bijection.

Let us consider the following two cases
\begin{enumerate}
\item[\rm (a)] $n$ is a positive odd  integer, $n > 1$,
\item[\rm (b)] $n =  2^{l}(\afrac {n}{2^{l}})$, where $l$ is a positive integer and $\afrac {n}{2^{l}}$ is odd.
\end{enumerate}

\bigskip
Case (a). Let $(k, m)$ be a fixed pair of positive odd integers such that $n = km$  and $m > 1$. We assume that $m > 1$, because $(n, 1, 0)$ is not a proper vector for $n > 1$. Hence, by condition (1) of Remark 4.2 and by condition (2) of Remark 4.1, there exists only one index-vector $(k, m, s)$ belonging to $\mathcal{M}_{n}$.

Since $l$ is odd, we see that $m$ is any divisor of $n$, $m > 1$. Hence, there are $\theta(n) - 1$ index-vectors belonging to $\mathcal{M}_{n}$.

\bigskip
Case (b). Let $(k, m)$ be a fixed pair of positive even integers such that $n = km$. Hence, by conditions  (1) and (2) of Remark 4.2 there exist only two index-vectors $(k, m, s_{1})$ and $(k, m , s_{2})$ belonging to $\mathcal{M}_{n}$, for some integers $s_{1}$ and $s_{2}$. 

Since $n =  2^{l}(\afrac {n}{2^{l}})$, $l > 0$ and $\afrac {n}{2^{l}}$ is odd, we see that $m = 2^{i}m_{j}$, where $1 \leqslant i  < l$ and $m_{j}$ is any divisor of $\afrac {n}{2^{l}}$. Hence, there are 
 $2(l-1) \theta (\afrac{n}{2^{l}})$  index-vectors belonging to  $\mathcal{M}_{n}$.
\bigskip

Let now $(k, m)$ be a fixed pair of positive integers such that $n = km$, $m$ is odd, $m > 1$ and $k$ is even. Hence, by condition (2) of Remark 4.2 and by condition (1) of Remark~4.1, there exists only one index-vector $(k, m, s)$ belonging to $\mathcal{M}_{n}$, for some integer $s$.

Since $n =  2^{l}(\afrac {n}{2^{l}})$, $l > 0$ and $\afrac {n}{2^{l}}$ is odd, we see that $m$ is a divisor of  $\afrac {n}{2^{l}}$, $m > 1$. Hence, there are $\theta (\afrac{n}{2^{l}}) - 1$ index-vectors  belonging to $\mathcal{M}_{n}$. 
\bigskip

Adding $2(l-1) \theta(\afrac{n}{2^{l}})$ to $\theta(\afrac{n}{2^{l}})-1$ we have $(2l - 1) \theta (\afrac{n}{2^{l}})-1$ index-vectors belonging to $\mathcal{M}_{n}$.  Hence, lemma holds because $\omega:\mathcal{S}_{n}  \rightarrow \mathcal{M}_{n}$ is a bijection.
\end{proof}
\bigskip

\section{Triangulations with codes of order $1$ or $2$}
\bigskip

We recall that $X_{n} = \{(k, m , s)\in \mathbb{\mathbb{Z}}^3: \ 1 \leqslant m \leqslant n, \  0 \leqslant s  < m \ \hbox{and} \ km = n\}$. 

In the following lemma, we state some conditions equivalent to the condition that triangulation $P$ has a code of order $2$.

\begin{lemma}\label{lemma6.0}
The following conditions are equivalent:
\begin{enumerate}
\item[\rm (1)] $P$ has an orbit of order $1$ which is not of the form $\{(\omega, \omega, 0)\}$ and $\{(\omega,  3\omega^{2}, \omega)\}$, for any integer $\omega$. 
\item[\rm (2)] $P$ has an orbit of order $1$ but $P$ does not have a code of order $1$.
\item[\rm (3)] $P$ has an orbit of order $1$ but $P$ is not symmetric,
\item[\rm (4)] $P$ has a code of order $2$.
\end{enumerate} 
\begin{proof}
$(1) \Rightarrow (2)$. Suppose, by contradiction, that $P$ has a code of order $1$. Then, by conditions $(2)-(3)$ of  Lemma~\ref{lemma5.2}, it has an orbit of the form $\{(\omega, \omega, 0)\}$ or $\{(\omega,  3\omega^{2}, \omega)\}$, for some integer $\omega$.

$(2) \Rightarrow (3)$ Suppose, by contradiction that $P$ is symmetric. Since $P$  has an orbit of order~$1$, $\bar{P}$ has an orbit of order $1$. Hence, $P$ has a code of order $1$.

$(3) \Rightarrow (4)$. If $P$  is not symmetric, then, by condition (2) of Proposition \ref{proposition1.4}, $P$ and~$\bar{P}$ have disjoint orbits. Since $P$ and $\bar{P}$ have orbits of the same order, they have orbits of order $1$. Thus, $P$ has a code of order $2$.

$(4) \Rightarrow (1)$. Suppose, by contradiction, that $P$ has an orbit of the form $\{(\omega, \omega, 0)\}$ or $\{(\omega,  3\omega^{2}, \omega)\}$, for some  integer $\omega$. Then, by conditions (2) and (3) of Lemma \ref{lemma5.2}, $P$ has a code of order $1$. 
\end{proof}
\end{lemma}

\begin{lemma}\label{lemma6.1}
Let $n = 2^{l}3^{\alpha}p^{\alpha_{1}}_{1}\ldots p^{\alpha_{u}}_{u}q^{\beta_{1}}_{1}\ldots q^{\beta_{w}}_{w}$ $(l, \alpha, \alpha_{i}, \beta_{i}\geqslant 0)$ be the decomposition of $n$ into primes such that $p_{i} \equiv\ 1 \pmod{ 3}$, for $i = 1,\ldots, u$, and $q_{i} \equiv 2 \pmod{ 3}$, for $i = 1,\ldots, w$. The following implications are true:
 \\[4pt] 
 If $n \neq {\omega}^{2}$ and $n  \neq 3{\omega}^{2}$ for any positive integer $\omega$, then there exist $\frac{\theta^{\star}(n)}{2}$  triangulations (non-isomorphic) in $\mathcal{P}_{n}$ with codes of order $2$ and no triangulations with code of order $1$.
  \\[4pt] 
If $n = {\omega}^{2}$ or $n = 3{\omega}^{2}$ for some positive integer $\omega$, then there exist $\frac{\theta^{\star}(n) - 1}{2}$    triangulations (non-isomorphic) in $\mathcal{P}_{n}$ with codes of order $2$ and one triangulation with code of order $1$,
 \\[4pt] 
 where
\[
\begin{array}{ll}
\theta^{\star} (n) = 
\left\{
\begin{array}{ll}
\theta(p^{\alpha_{1}}_{1}\ldots p^{\alpha_{u}}_{u}), & \hbox{if  } \ \l \ \hbox{and} \ \beta_{i}\ \hbox{are even},  \ \hbox{for every} \ i = 1, \ldots,  w, 
\\
 0, &\hbox{if}\  l \ \hbox{or} \ \beta_{i} \ \hbox{is odd}, \ \hbox{for some}\ i = 1, \ldots,  w.
\end{array}
\right.
\end{array}
\]
\end{lemma}
\begin{proof}  Let $n = 2^{l}3^{\alpha}p^{\alpha_{1}}_{1}\ldots p^{\alpha_{u}}_{u}q^{\beta_{1}}_{1}\ldots q^{\beta_{w}}_{w}$ $(l, \alpha, \alpha_{i}, \beta_{i}\geqslant 0)$ be the decomposition of $n$ into primes such that $p_{i} \equiv\ 1 \pmod{ 3}$, for $i = 1,\ldots, u$, and $q_{i} \equiv 2 \pmod{ 3}$, for $i = 1,\ldots, w$.

Let $t(n)$, $n \in N$, be the number of solutions of the congruence $$x^{2} + x + 1  \equiv 0 \pmod{n}\ \hbox{for} \ 0 \leqslant x < n.$$ 
\noindent
Certainly:  
$t(1) = 1$, $t( 2^{l}) = 0$ (for $ l > 0$), $t(3) = 1$ and $t(3^{\alpha}) =0$ (for $\alpha > 1$).
 
\noindent
From Mohar theorem \cite[Proposition 3.2]{flo9} (see Proposition \ref{prop1.7}) it follows that 

$t(p^{s}) = 2$, where $p$ is a prime, $p \equiv\ 1 \pmod{ 3}$ and $s > 0$,

$t(q^{s}) = 0$, where $q$ is a prime, $q \equiv\ 2  \pmod{ 3}$ and $s > 0$ and
\begin{align}
t(2^{l}3^{\alpha}p^{\alpha_{1}}_{1}\ldots p^{\alpha_{u}}_{u}q^{\beta_{1}}_{1}\ldots q^{\beta_{w}}_{w}) = t(2^{l})t(3^{\alpha})t(p^{\alpha_{1}}_{1})\ldots t(p^{\alpha_{u}}_{u})t(q^{\beta_{1}}_{1})\ldots t(q^{\beta_{w}}_{w}).
\end{align}

From Proposition \ref{prop1.6} we conclude that:  $\{(k, m, s)\} \in X_{n}$ is the orbit of some triangulation $P \in \mathcal{P}_{n}$ if and only if  $(k, m, s) = (k, kz, kx)$, where $k^{2}z = n$ and $0 \leqslant x < \afrac{n}{k^{2}}$ is a solution of the congruence $x^{2} + x + 1 \equiv 0 \pmod{\afrac{n}{k^{2}}}$. Hence, 
$$\sum\limits_{k^{2}|n} t(\frac {n}{k^{2}}) \  \hbox{ is the number of orbits of order} \ 1 \ \hbox{of  triangulations in} \ \mathcal{P}_{n} ,$$
where $\sum$ is taken over all divisors $k^{2}$ of  $n$. 

Let us consider the case that $ l$ or $\beta_{i}$ is odd, for some $i = 1, \ldots,  w$. By Proposition~\ref{prop1.7}, $ t(\afrac {n}{k^{2}}) = 0$, for every divisor $k^{2}$ of $n$. Thus, there exists no triangulation of $\mathcal{P}_{n}$ with orbit of order $1$. Since $\theta^{\star}(n) = 0$, Lemma \ref{lemma6.1} holds.

Now let us  consider the case where $ l$ and $\beta_{i}$ is even, for every $i = 1, \ldots,  w$. Let us denote $m =p^{\alpha_{1}}_{1}\ldots p^{\alpha_{u}}$. Assume that $\alpha_i$ is odd (even) for $1 \leqslant i \leqslant r$ ($r+1 \leqslant i \leqslant u$, respectively). Hence, by (4), we obtain
\begin{align*}
&\sum\limits_{k^{2}|n} t(\frac {n}{k^{2}})  =  \sum\limits_{k^{2}|m} t(\frac {p^{\alpha_{1}}_{1}\ldots p^{\alpha_{u}}_{u}}{k^{2}}) =\sum\limits_{k^{2}|m} t(\frac {p^{\alpha_{1}}_{1}\ldots p^{\alpha_{r}}_{r}p^{\alpha_{r+1}}_{r+1}\ldots p^{\alpha_{u}}_{u}}{k^{2}}) =
 \\[4pt]
& 
 \sum t(p^{s_{1}}_{1}\ldots p^{s_{r}}_{r} p^{s_{r+1}}_{r+1} \ldots p^{s_{u}}_{u}) ,\ \hbox{ where } \ \sum \hbox{ is taken over all}
\end{align*} 
\\[1pt]
$1 \leqslant s_{i} \ odd  \leqslant \alpha_{i}$,  for $1\leqslant i \leqslant r$, and $0 \leqslant s_{i} \  even  \leqslant \alpha_{i}$, for $r+1 \leqslant i  \leqslant u$. Hence, by (4), we obtain

\begin{align*}
&\sum\limits_{k^{2}|n} t(\frac {n}{k^{2}})  = 
\\[4pt]
& \sum\limits_{1 \leqslant s_{1} odd  \leqslant \alpha_{1}}t(p^{s_{1}}_{1}) \ldots  
\sum\limits_{1 \leqslant s_{r} odd  \leqslant \alpha_{r}} t(p^{s_{r}}_{r})
 \ \sum\limits_{0  \leqslant s_{r+1} even  \leqslant \alpha_{r+1}}t(p^{s_{r+1}}_{r+1})\ldots 
 \sum\limits_{0  \leqslant s_{u} even  \leqslant \alpha_{u}} t( p^{s_{u}}_{u}) =
 \\[4pt]
& 2 \frac{(\alpha_{1}+1)}{2} \ldots 2\frac{(\alpha_{r}+1)}{2}\{2\frac{\alpha_{r+1}}{2} + 1\} \ldots \{2\frac{\alpha_{u}}{2}+1\} =  \theta(p^{\alpha_{1}}_{1}\ldots p^{\alpha_{u}}_{u})
\end{align*}
Certainly, if $\alpha_i$ is odd (even) for every $1 \leqslant i \leqslant u$, then
we obtain the same equality as above.

Assume that $n \neq {\omega}^{2}$ and $n  \neq 3{\omega}^{2}$, for any positive integer $\omega$. If $P \in\mathcal{P}_{n}$ has an orbit of order $1$, then, it is not of the form $\{(\omega, \omega, 0)\}$ or $\{(\omega, 3{\omega}, \omega)\}$. Thus, by conditions (1) and (4) of Lemma \ref{lemma6.0}, $P$ has a code of order $2$. Hence, we see that
$$\frac {1}{2}\sum\limits_{k^{2}|n} t(\frac {n}{k^{2}}) = \frac{\theta^{\star}(n)}{2} $$
 is the number of triangulations in  $\mathcal{P}_{n}$ with codes of order $2$.

Let now $n = {\omega}^{2}$ (or $n  = 3{\omega}^{2}$), for some positive integer $\omega$. Then, by conditions $(2)-(3)$ of  Lemma \ref{lemma5.2}, a triangulation with an orbit of the form $\{(\omega, \omega, 0)\}$ or $\{(\omega, 3{\omega}, \omega)\}$, for some integer $\omega$,  is the only one triangulation in $\mathcal{P}_{n}$ with a code of order $1$. If $P \in\mathcal{P}_{n}$ has an orbit of order $1$ different from $\{(\omega, \omega, 0)\}$ (or $\{(\omega, 3{\omega}, \omega)\}$, respectively), then, by conditions (1) and (4) of Lemma \ref{lemma6.0},  $P$ has a code of order $2$. Hence, we obtain
$$\frac {1}{2}\sum\limits_{k^{2}|n, \ k^{2} \neq n} t(\frac {n}{k^{2}}) = \frac {1}{2}\{\sum\limits_{k^{2}|n} t(\frac {n}{k^{2}}) - t(1)\} = \frac {\theta^{\star}(n) - 1}{2}$$
$$(\hbox{or}  \ \frac {1}{2}\sum\limits_{k^{2}|n, \ 3k^{2} \neq n} t(\frac {n}{k^{2}}) = \frac {1}{2}\{\sum\limits_{k^{2}|n} t(\frac {n}{k^{2}} - t(3)\} = \frac{\theta^{\star}(n) - 1}{2})$$
is the number of triangulations in $\mathcal{P}_{n}$ with codes of order $2$.
\end{proof}
\bigskip

\section{The enumeration of triangulations}
\bigskip

Notice that if $(k, m, s) \in X_{n}$, then $m$ is a divisor of $n$, $0 \leqslant s < m$ and $k = \afrac{n}{m}$. Hence  $|X_{n}| = \sigma(n)$. Since, by Remark \ref{rem4.1}, each vector of~$X_{n}$  different from $(n,1,0)$, $(1,n, n-1)$ and $(1,n, 0)$ (for $n > 1$) is an index-vector of some triangulation in $\mathcal{P}_{n}$, there are $\sigma(n) - 3$ index-vectors of triangulations in $\mathcal{P}_n$ (for $n > 1$). 

In the following lemma, we state some conditions equivalent to the condition that triangulation $P$ has a code of order $3$.
 
\begin{lemma}\label{lemma7.1}
The following conditions are equivalent:
\begin{enumerate}
\item[\rm (1)] $P$ is symmetric but it does not have orbit of the form $\{(\omega, \omega, 0)\}$ and $\{(\omega,  3\omega^{2}, \omega)\}$, for any integer $\omega$.
\item[\rm (2)] $P$ is symmetric and it does not have a code of order $1$,
\item[\rm (3)] $P$ is symmetric and it does not have an orbit of order $1$,
\item[\rm (4)] $P$ has a code of order $3$.
\end{enumerate} 
\begin{proof}
$(1) \Rightarrow (2)$ Suppose, by contradiction, that $P$ has a code of order $1$. Then, by conditions $(2)$ and $(3)$ of Lemma \ref{lemma5.2}, it has an orbit of the form $\{(\omega, \omega, 0)\}$ or $\{(\omega,  3\omega^{2}, \omega)\}$, for some integer $\omega$.

$(2) \Rightarrow (3)$. Suppose, by contradiction, that $P$ has an orbit of order $1$. Then, $\bar{P}$ has an orbit of order $1$. Since, $P$ is symmetric, $P$ has a code of order $1$.

$(3) \Rightarrow (4)$. If $P$ has no orbit of order $1$, then, by condition $(3)$ of Proposition 1.4, $P$ has an orbit of order $3$.
 Since, $P$ is symmetric, $P$ has a code of order $3$.

$(4) \Rightarrow (1)$. Suppose, by contradiction, that $P$ has an orbit of the form form $\{(\omega, \omega, 0)\}$ or $\{(\omega,  3\omega^{2}, \omega)\}$, for some integer $\omega$. Then. by conditions $(2)$ and $(3)$ of Lemma \ref{lemma5.2}, $P$ has a code of order $1$.
\end{proof}
\end{lemma}
\bigskip

\begin{proof} [\bf Proof of Theorem \ref{theo1.1}]
Let $d(n)$ be the number of all (non-isomorphic) triangulations in  $\mathcal{P}_{n}$. Notice that, by condition (1) of Proposition \ref{prop1.5}, $d(n)$ is also the number of different codes of triangulations in $\mathcal{P}_{n}$. 

Let $n = 2^{l}3^{\alpha}p^{\alpha_{1}}_{1}\ldots p^{\alpha_{u}}_{u}q^{\beta_{1}}_{1}\ldots q^{\beta_{w}}_{w}$, $n > 1$ $(l, \alpha, \alpha_{i}, \beta_{i}\geqslant 0)$, be the decomposition of $n$ into primes such that $p_{i} \equiv 1 \pmod{ 3}$, for $i = 1,\ldots, u$, and $q_{i} \equiv 2 \pmod{ 3}$, for $i = 1,\ldots, w$.  Let us consider the following cases:
\begin{enumerate}
\item[\rm (a)] $l = 0$ and $n \neq 3\omega^{2}$ and $n \neq \omega^{2}$, for any integer $\omega$,
\item[\rm (b)]  $l = 0$, and $n = 3\omega^{2}$ or $n = \omega^{2}$, for some  integer $\omega$,
\item[\rm (c)]  $l > 0$ and $n \neq 3\omega^{2}$ and $n \neq 3\omega^{2}$, for any  integer $\omega$,
\item[\rm (d)] $l > 0$, and $n = 3\omega^{2}$ or $n = \omega^{2}$, for some integer $\omega$.
\end{enumerate} 
\bigskip

Case (a). 
By Lemma \ref{lemma6.1}, there exist $\afrac{\theta^{\star}(n)}{2}$ codes of order $2$ of triangulations in $\mathcal{P}_n$, and there is no code of order $1$. 

Since $n$ is odd, by Lemma \ref{lemma5.3}, there exist $\theta(n) - 1$ symmetric triangulations in $\mathcal{P}_n$. Since no one of them has a code of order $1$, by conditions (2) and (4) of Lemma \ref{lemma7.1},  there exist $\theta(n) - 1$ codes of order $3$. Hence, by condition (2)-(3) of Proposition \ref{prop1.5}, there are 
\begin{align*}
d(n) &= \frac{\sigma(n) - 3 -  3(\theta(n)- 1) - \theta^{\star}(n)}{6} + \theta(n) - 1  + \frac{\theta^{\star}(n)}{2} 
\\[4pt]
 &= \frac{\sigma(n) +  3\theta(n) + 2 \theta^{\star}(n)}{6} -1
\end{align*} 
  codes of triangulations in $\mathcal{P}_n$. Thus, condition  (2) of Theorem \ref{theo1.1} holds.
\bigskip

Case (b). By Lemma \ref{lemma6.1}, there exist $\afrac {(\theta^{\star}(n) - 1)}{2}$ codes of order $2$  of  triangulations in $\mathcal{P}_n$ and one code of order $1$ ($\{(\omega, \omega, 0)\}$ for $\alpha$ even, or $\{(\omega, 3\omega, \omega)\}$ for $\alpha$ odd).  

Since $n$ is odd, by Lemma \ref{lemma5.3}, there exist $\theta(n) - 1$ symmetric triangulations in $\mathcal{P}_n$.  Since only one of them has a code of order $1$,  by conditions (2) and (4) of Lemma \ref{lemma7.1}, there exist  $\theta(n) - 2$ codes of order~$3$. Hence, by conditions (2)-(3) of Proposition \ref{prop1.5}, there are
\begin{align*}
d(n) &= \frac{\sigma(n) - 3 -  3(\theta(n) -2) - 1 -  (\theta^{\star}(n) -1)}{6} + \theta(n)  - 1 + \frac{\theta^{\star}(n) - 1}{2}
\\[4pt]
&=  \frac{\sigma(n) +  3\theta(n) + 2 \theta^{\star}(n)}{6} -1 
\end{align*}
 codes of triangulations in $\mathcal{P}_n$. Thus, condition (2) of Theorem \ref{theo1.1} holds.
\bigskip

Case (c). By Lemma \ref{lemma6.1}, there exist $\afrac{\theta^{\star}(n)}{2}$ codes of order $2$ of triangulations in $\mathcal{P}_n$, and there is no code of order $1$.
 
 Since $n = 2^{l}(\afrac {n}{2^{l}})$, $l > 0$ and $\afrac {n}{2^{l}}$ is odd, by Lemma \ref{lemma5.3}, there exist $(2l - 1) \theta (\afrac{n}{2^{l}}) - 1$ symmetric triangulations in $\mathcal{P}_n$. Therefore,  by conditions (2) and (4) of Lemma \ref{lemma7.1}, there exist $(2l - 1)\theta(\afrac{n}{2^{l}}) - 1$  codes of order~$3$.
Hence, by conditions (2)-(3) of Proposition \ref{prop1.5}, there are 

\begin{align*} 
d(n) &= \frac{\sigma(n) - 3 - 3[(2l - 1) \theta (\afrac{n}{2^{l}})- 1] - \theta^{\star}(n)}{6} +  (2l - 1)\theta (\afrac{n}{2^{l}}) -1+\frac {\theta^{\star}(n)}{2}
\\[4pt]
&=   \frac{\sigma(n) + 3(2l - 1)\theta(\afrac{n}{2^{l}}) + 2 \theta^{\star}(n)}{6} - 1
\end{align*}
 codes  of triangulations in $\mathcal{P}_n$. Hence, condition (3)  of Theorem \ref{theo1.1} holds.
\bigskip

Case (d). By Lemma \ref{lemma6.1}, there exist $\afrac {(\theta^{\star}(n) - 1)}{2}$  codes of order $2$  of  triangulations in $\mathcal{P}_n$ and one code of order $1$ ($\{(\omega, \omega, 0)\}$ for $\alpha$ even, or $\{(\omega, 3\omega, \omega)\}$ for $\alpha$ odd).  

 Since $n = 2^{l}(\afrac {n}{2^{l}})$, $l > 0$ and $\afrac {n}{2^{l}}$ is odd, by Lemma \ref{lemma5.3}, there are ${(2l - 1) \theta (\afrac{n}{2^{l}}) - 1}$ symmetric triangulations in $\mathcal{P}_n$. Since only one of them has a code of order $1$, by conditions (2), (4) of Lemma \ref{lemma7.1}, there exist $(2l - 1) \theta (\afrac{n}{2^{l}}) - 2$ codes of order  $3$. Hence, by conditions (2)-(3) of Proposition \ref{prop1.5}, there are
\begin{align*} 
d(n) &= \frac{\sigma(n)- 3 -  3[(2l - 1) \theta (\afrac{n}{2^{l}}) - 2] - 1 - (\theta^{\star}(n) -1)}{6} + (2l - 1)\theta(\afrac{n}{2^{l}}) - 1 
\\[4pt]
&+ \frac {\theta^{\star}(n) - 1}{2}=  \frac{\sigma(n) + 3(2l - 1)\theta(\afrac{n}{2^{l}}) + 2 \theta^{\star}(n)}{6} - 1
\end{align*}
 codes  of triangulations in $\mathcal{P}_n$. Thus, condition (3)  of Theorem \ref{theo1.1} holds.
\end{proof}

\section{Odd perfect numbers}
\bigskip

\noindent
Let  $n$ be an odd positive integer and suppose that $n = p^{\alpha_{1}}_{1}\ldots p^{\alpha_{u}}_{u}q^{\beta_{1}}_{1}\ldots q^{\beta_{w}}_{w}$ $(\alpha_{i}, \beta_{j}\geqslant~0)$ is the decomposition of $n$ into primes such that $p_{i} \equiv 1 \pmod{ 3}$, for $i = 1,\ldots, u$, and $q_{j} \equiv 2 \pmod{ 3}$, for $j = 1,\ldots, w$. By Theorem \ref{theo1.1} we have the following equality: 
\bigskip

\begin{align} 
\sigma(n) +  3\theta(n) + 2 \theta^{\star}(n) \equiv 0 \pmod{6},
\end{align}
where
\[
\begin{array}{ll}
\theta^{\star} (n) = 
\left\{
\begin{array}{ll}
\theta(p^{\alpha_{1}}_{1}\ldots p^{\alpha_{u}}_{u}), &\hbox{if} \  \ \beta_{j} \ \hbox{are even}, \ \hbox{for every} \ j = 1, \ldots,  w, 
\\
 0, &\hbox{if}\ \ \beta_{j} \ \hbox{is odd}, \ \hbox{for some}\ j = 1, \ldots,  w.
\end{array}
\right.
\end{array}
\]
\bigskip
 
It follows from the Euler theorem (see Theorem 1.3) that, if $n$ is an odd perfect number, then
\begin{align} 
n = \pi^{\alpha}M^{2},\hbox { where }
\pi \hbox{ is prime }, gcd(\pi,M) = 1 \hbox{  and } \pi  \equiv  \alpha  \equiv  n \equiv 1 \pmod{4}.
\end{align}

The following lemma is well known (see  Starni \cite{flo14}). We prove it 
using equality (5). 
\begin{lemma} \label{lemma7.1} Let $n$ be an odd perfect number number satisfying condition $(6)$. The following conditions are equivalent. 
\begin{enumerate}
\item [$(i)$] $3 \nmid n$,
\item [$(ii)$] $\pi  \equiv 1 \pmod{3}$,
\item [$(iii)$] $\pi  \equiv 1 \pmod{12}$,
\item [$(iv)$] $n \equiv 1 \pmod{12}$.
\end{enumerate}
\end{lemma}
\begin{proof} $(i) \Rightarrow (ii)$. Assume that $\pi  \equiv 2 \pmod{3}$. Hence,  $\pi = q_{j_{0}}$ and $\alpha = \beta_{j_{0}}$, for some $1 \leqslant j_{0} \leqslant w$. By condition (6), $\alpha = \beta_{j_{0}}$ is odd. Hence, $\theta^{\star} (n) =0$. Moreover,  $\theta(n)$ is even, because $\beta_{j_{0}}$ is odd. Therefore, by equqlity (5), we have
\[2n = \sigma(n)  \equiv 0 \pmod{6}.\] 

$(ii) \Rightarrow (iii)$. By  (6), $\pi  \equiv 1 \pmod{4}$. If $\pi  \equiv 1 \pmod{3}$, then $\pi  \equiv 1 \pmod{12}$.

$(iii) \Rightarrow (iv)$. Notice that $M^{2} \equiv 1  \pmod{3}$ and $M^{2} \equiv 1  \pmod{4}$. Hence,  $M^{2} \equiv 1  \pmod{12}$. If $\pi \equiv 1  \pmod{12}$, then $n = \pi^{\alpha}M^{2} \equiv 1  \pmod{12}$.

$(iv) \Rightarrow (i)$. Obvious.
\end{proof}

The following proposition was proved by Touchard \cite{flo13}. It follows from the Euler theorem and Lemma \ref{lemma7.1}.

\begin{proposition} \label{prop.7.1} $[Touchard]$ If $n$ is an odd perfect number, then \[n  \equiv 9\pmod{36}\ \hbox{ or }\  n  \equiv 1 \pmod{12}.\]
\end{proposition} 

\begin{proof} [\bf Proof of Theorem \ref{theo1.4}]

 Let $n$ be an odd perfect number satisfying condition (6) and suppose that $3 \nmid n$. Hence,  by Lemma~\ref{lemma7.1}, $\pi \equiv n \equiv 1  \pmod{3}$. Hence, $\pi = p_{i_{0}}$ and $\alpha = \alpha_{i_{0}}$, for some $1 \leqslant i_{0} \leqslant u$. By condition (6), $\beta_{i}$ are even, for every $i = 1, \ldots,  w$, and $\alpha$ is odd.  Hence,
\[\theta^{\star} (n)  = \theta(p^{\alpha_{1}}_{1}\ldots p^{\alpha_{u}}_{u}) \ \hbox{ and }\ \theta(n) \ \hbox{is even}.\]
Thus,  by equality (5), 

\[2n + 2  \theta(p^{\alpha_{1}}_{1}\ldots p^{\alpha_{u}}_{u}) \equiv 0 \pmod{6}.\]

Since $n \equiv 1  \pmod{3}$, 
\begin{align} 
\theta(p^{\alpha_{1}}_{1}\ldots0p^{\alpha_{u}}_{u})  \equiv 2 \pmod{3}.
\end{align} 
By condition (6), $\alpha_{i_{0}} = \alpha = 4k +1$, for some integer $k$. Moreover, $\alpha_{i}$ is even, for $i \neq i_{0}$ and $i = 1,\ldots, u$. Hence, it is sufficient to consider the following two cases:
\begin{enumerate}
\item [$(i)$] $\theta(p^{\alpha_{1}}_{1}\ldots p^{\alpha_{u}}_{u}) = (4k +2) (3l+1)$, for some integer $l$,
\item  [$(ii)$] $\theta(p^{\alpha_{1}}_{1}\ldots p^{\alpha_{u}}_{u}) = (4k +2) (3l+2)$, for some integer $l$,
\end{enumerate}
\bigskip

Case (i). $\theta(p^{\alpha_{1}}_{1}\ldots p^{\alpha_{u}}_{u}) = (4k +2) (3l+1) \equiv 4k +2 \pmod{3}$. Hence, by equality (7), $k = 3s$ and $\alpha = 12s + 1$, for some integer $s$. 
\bigskip

Case (ii). $\theta(p^{\alpha_{1}}_{1}\ldots p^{\alpha_{u}}_{u}) = (4k +2) (3l+2) \equiv 8k +4 \pmod{3}$. Hence, by equality (7), $k = 3s +2$ and $\alpha = 4(3s +2) + 1 = 12s +9$, for some integer $s$. 
\bigskip

Finally, by Lemma 7.1, $\pi  \equiv  n \equiv 1 \pmod{12}$.
\end{proof}



\end{document}